\documentclass[12pt]{article} 
\usepackage{amssymb,amsmath}
\usepackage{amsthm}
\input{psfig.sty}

\newtheorem{theorem}{Theorem}[section]

\newtheorem{lemma}[theorem]{Lemma}
\newtheorem{corollary}[theorem]{Corollary}

\theoremstyle{definition}

\theoremstyle{remark}

\begin{document}
\newcommand{\beq}{\begin{equation}} \newcommand{\eeq}{\end{equation}}
\newcommand{\zz}{\mathbb{Z}}
\newcommand{\pp}{\mathbb{P}} 
\newcommand{\nn}{\mathbb{N}}
\newcommand{\rr}{\mathbb{R}}
\newcommand{\bm}[1]{{\mbox{\boldmath $#1$}}}
\newcommand{\con}{\mathrm{Comp}(n)}
\newcommand{\sn}{\mathfrak{S}_n} 
\newcommand{\fs}{\mathfrak{S}}
\newcommand{\st}{\,:\,} 
\newcommand{\as}{\mathrm{as}}
\newcommand{\is}{\mathrm{is}}
\newcommand{\lgn}{\mathrm{len}}
\newcommand{\dis}{\displaystyle}

\thispagestyle{empty}

\vskip 20pt
\begin{center}
{\large\bf Longest Alternating Subsequences of
  Permutations}
\vskip 15pt
{\bf Richard P. Stanley}\\
{\it Department of Mathematics, Massachusetts Institute of
Technology}\\
{\it Cambridge, MA 02139, USA}\\
{\texttt{rstan@math.mit.edu}}\\[.2in]
{\bf\small version of 15 November 2005}\\
\end{center}

\begin{abstract}
The length $\is(w)$ of the longest increasing subsequence of a
permutation $w$ in the symmetric group $\sn$ has been the object of
much investigation. We develop comparable results for the length
$\as(w)$ of the longest alternating subsequence of $w$, where a
sequence $a,b,c,d,\dots$ is \emph{alternating} if
$a>b<c>d<\cdots$. For instance, the expected value (mean) of $\as(w)$
for $w\in\sn$ is exactly $(4n+1)/6$ if $n\geq 2$.
\end{abstract}

\section{Introduction.} \label{sec1}
\indent Let $\sn$ denote the symmetric group of permutations of
$1,2,\dots,n$, and let $w=w_1\cdots w_n\in\sn$. An \emph{increasing
  subsequence} of $w$ of length $k$ is a subsequence $w_{i_1}\cdots
w_{i_k}$ satisfying
  \[ w_{i_1} < w_{i_2}<\cdots <w_{i_k}. \]
There has been much recent work on the length is$_n(w)$
of the longest increasing subsequence of a permutation $w\in\sn$. A
highlight is the asymptotic determination of the expectation $E(n)$ of
$\is_n$ by Logan-Shepp \cite{l-s} and Vershik-Kerov \cite{v-k}, viz.,
 \beq E(n):= \frac{1}{n!}\sum_{w\in\sn} \is_n(w)\sim 2\sqrt{n},\ 
    n\rightarrow \infty. \label{eq:en} \eeq
Baik, Deift and Johansson \cite{b-d-j} obtained a vast strengthening
of this result, viz., the limiting distribution of is$_n(w)$ as
$n\rightarrow \infty$. Namely, for $w$ chosen uniformly from $\sn$ we
have  
  \beq \lim_{n\rightarrow\infty}
  \mathrm{Prob}\left(\frac{\is_n(w)-2\sqrt{n}}{n^{1/6}}\leq t\right)
     = F(t), \label{eq:bdj} \eeq
where $F(t)$ is the Tracy-Widom distribution. The proof uses a result
of Gessel \cite{gessel} that gives a generating function for the
quantity 
  \[ u_k(n) = \#\{w\in\sn\st \is(w)\leq k\}. \]
Namely, define
   \begin{align*}   U_k(x) & = \sum_{n\geq 0} u_k(n)\frac{x^{2n}}{n!^2},\
   k\ge 1\\  
    I_i(2x) & = \sum_{n\geq 0}\frac{x^{2n+i}}{n!\,(n+i)!},\
       i\in\zz. \end{align*}
The function $I_i$ is the \emph{hyperbolic Bessel function} of the
first kind of order $i$. Note that $I_i(2x)=I_{-i}(2x)$. Gessel then
showed that
  \[ U_k(x) = \det\left( I_{i-j}(2x)\right)_{i,j=1}^k. \]
\indent In this paper we will develop an analogous theory for
\emph{alternating subsequences}, i.e., subsequences $w_{i_1}\cdots
w_{i_k}$ of $w$ satisfying
  \[ w_{i_1}>w_{i_2}<w_{i_3}>w_{i_4}<\cdots\,w_{i_k}. \]
Note that according to our definition, an alternating sequence
$a,b,c,\dots$ (of length at least two) must begin with a descent
$a>b$. Let $\as(w)=\as_n(w)$ denote the length (number of terms) of the
longest alternating subsequence of $w\in\sn$, and let
  \[ a_k(n) =\#\{w\in\sn\st \as(w)=k\}. \]
For instance, $a_1(w)=1$, corresponding to the permutation $12\cdots
n$, while $a_n(n)$ is the total number of alternating permutations in
$\sn$. This number is customarily denoted $E_n$. A celebrated result
of Andr\'e \cite{andre}\cite[{\S}3.16]{ec2} states that
  \beq \sum_{n\geq 0} E_n\frac{x^n}{n!} = \sec x + \tan x. 
   \label{eq:engf} \eeq
The numbers $E_n$ were first considered by Euler (using (\ref{eq:engf})
as their definition) and are known as \emph{Euler numbers}. Because of
(\ref{eq:engf}) $E_{2n}$ is also known as a \emph{secant number} and
$E_{2n-1}$ as a \emph{tangent number}.  

Define 
  \begin{align} b_k(n) & = \#\{w\in\sn\st \as(w)\leq k\} \nonumber\\
   & = a_1(n) + a_2(n) + \cdots + a_k(n), \label{eq:bkak} \end{align}
so for instance $b_k(n)=n!$ for $k\geq n$. Also define the generating
functions
  \begin{align} A(x,t) & = \sum_{k,n\geq 0}a_k(n)t^k\frac{x^n}{n!}
  \label{eq:axt} \\
   B(x,t) & = \sum_{k,n\geq 0}b_k(n)t^k\frac{x^n}{n!}. \nonumber
   \end{align} 
Our main result (Theorem~\ref{thm:main}) is the formulas
    \begin{align} B(x,t) & = \frac{1+\rho+2te^{\rho x}+(1-\rho)
     e^{2\rho x}}{1+\rho-t^2+(1-\rho-t^2)e^{2\rho x}}
    \label{eq:fxt}\\[.05in] A(x,t) & =(1-t)B(x,t), \nonumber \end{align}
where $\rho=\sqrt{1-t^2}$. 

As a consequence of these formulas we obtain explicit formulas for
$a_k(n)$ and $b_k(n)$:
  \begin{align*} 
  b_k(n) & = \frac{1}{2^{k-1}} \sum_{\substack{r+2s\leq k\\ r\equiv
     k\,(\mathrm{mod}\,2)}} 
     (-2)^s{k-s\choose (k+r)/2}{n\choose s}r^n\\ 
   a_k(n) & = b_k(n)-b_{k-1}(n).\end{align*}
We also obtain from equation (\ref{eq:fxt}) formulas for the factorial
moments
  \[ \nu_k(n) = \frac{1}{n!} \sum_{w\in\sn}\as(w)(\as(w)-1)\cdots
     (\as(w)-k+1). \]
For instance, the mean $\nu_1(n)$ and variance var$(\as_n)=
\nu_2(n)+\nu_1(n)-\nu_1(n)^2$ are given by 
   \begin{align} \label{eq:mvar} \begin{split} \nu_1(n) & =
     \dis\frac{4n+1}{6},\ n>1\\[.05in] 
        \mathrm{var}(\as_n) & =  \dis\frac{8}{45}n-\frac{13}{180},\
        n\geq 4. \end{split}\end{align}
\indent The limiting distribution of $\as_n$ (the analogue of
equation~\ref{eq:bdj})) was obtained independently by Pemantle and
Widom, as discussed at the end of Section~\ref{sec3}. Rather than the 
Tracy-Widom distribution as in (\ref{eq:bdj}), this time we obtain a
Gaussian distribution.

\textsc{Note.} We can give an alternative description of $b_k(n)$
in terms of pattern avoidance. If $v=v_1 v_2\cdots v_k\in\fs_k$,
then we say that a permutation $w=w_1 w_2\cdots w_n\in\sn$
\emph{avoids} $v$ if $w$ has no subsequence $w_{i_1}w_{i_2}\cdots
w_{i_k}$ whose terms are in the same relative order as $v$
\cite[Ch.~4.5]{bona}\cite[{\S}7]{ids}. If $X\subset \fs_k$, then we
say that $w\in\sn$ \emph{avoids} $X$ if $w$ avoids all $v\in X$. Now
note that $b_{k-1}(n)$ is the number of permutations $w\in \sn$ that
avoid all $E_k$ alternating permutations in $\fs_k$.

After seeing the first draft of this paper Mikl\'os B\'ona pointed out
that the statistic $\as_n$ can be expressed very simply in terms of a
previously considered statistic on $\sn$, viz., the number of
\emph{alternating runs}. Hence our results can also be deduced from
known results on alternating runs. This development is discussed
further in Section~\ref{sec4}. In particular, it follows from
\cite{wilf} that the polynomials $T_n(t)=\sum_k a_k(n)t^k$ have
interlacing real zeros. This result can be used to give a third proof
(in addition to the proofs of Pemantle and Widom) that the limiting
distribution of $\as_n$ is Gaussian.

\section{The main generating function.}
The key result that allows us to obtain explicit formulas is the
following lemma.

\begin{lemma} \label{lemma:key}
Let $w\in\sn$. Then there is an alternating subsequence of $w$ of
maximum length that contains $n$.
\end{lemma}

\proof Let $a_1>a_2<\cdots a_k$ be an alternating subsequence of $w$
of maximum length $k=\as(w)$, and suppose that $n$ is not a term of
this subsequence. If $n$ precedes $a_1$ in $w$, then we can replace
$a_1$ by $n$ and obtain an alternating subsequence of length $k$
containing $n$. If $n$ appears between $a_i$ and $a_{i+1}$ in $w$,
then we can similarly replace the larger of $a_i$ and $a_{i+1}$ by
$n$. Finally, suppose that $n$ appears to the right of $a_k$. If $k$
is even that we can append $n$ to the end of the subsequence to obtain
a longer alternating subsequence, contradicting the definition of
$k$. But if $k$ is odd, then we can replace $a_k$ by $n$, again
obtaining an alternating subsequence of length $k$ containing
$n$. \qed

\medskip
We can use Lemma~\ref{lemma:key} to obtain a recurrence for
$a_k(n)$, beginning with the initial condition $a_0(0)=1$.

\begin{lemma} \label{lemma:rec}
Let $1\leq k\leq n+1$. Then
  \beq a_k(n+1) = \sum_{j=0}^{n} {n\choose j}
   \sum_{\substack{2r+s=k-1\\ r,s\geq
   0}}(a_{2r}(j)+a_{2r+1}(j))a_s(n-j).  
   \label{eq:aknrec} \eeq
\end{lemma}

\proof We can choose a permutation $w=a_1\cdots a_{n+1}\in\fs_{n+1}$
such that $\as(w)=k$ as follows. First choose $0\leq j\leq n$ such
that $a_{j+1}=n+1$. Then choose in ${n\choose j}$ ways the set
$\{a_1,\dots, a_j\}$. For $s\geq 0$ we can choose in $a_s(n-j)$
ways a permutation $w'=a_{j+2}\cdots a_{n+1}$ satisfying
$\as(w')=s$. Next we choose a permutation $w''=a_1\cdots a_j$ such
that the longest \emph{even} length of an alternating subsequence of
$w''$ is $2r=k-1-s$. We can choose $w''$ to satisfy
either $\as(w'')=2r$ or $\as(w'')=2r+1$. The concatenation
$w=w''(n+1)w'\in\fs_{n+1}$ will then satisfy $\as(w)=k$, and
conversely all such $w$ arise in this way. Hence equation
(\ref{eq:aknrec}) follows. 
\qed

\medskip
Now write
  \[ F_k(x) = \sum_{n\geq 0}a_k(n)\frac{x^n}{n!}. \]
For instance, $F_0(x)=1$ and $F_1(x)=e^x-1$.
Multiplying (\ref{eq:aknrec}) by $x^n/n!$ and summing on $n\geq 0$
gives 
  \beq F'_k(x) = \sum_{2r+s=k-1} (F_{2r}(x)+F_{2r+1}(x))F_s(x). 
    \label{eq:akx} \eeq
Note that 
  \[ A(x,t) = \sum_{k\geq 0} F_k(x)t^k, \]
where $A(x,t)$ is defined by (\ref{eq:axt}). Since $k-1-s$ is even in
(\ref{eq:akx}), we need to work with the even part $A_e(x,t)$ and odd
part $A_o(x,t)$ of $A(x,t)$, defined by
  \begin{equation} \label{eq:aeao}
   \begin{split} A_e(x,t) & =  \sum_{k\geq 0} F_{2k}(x)t^{2k}\\ 
     & =  \frac 12(A(x,t)+A(x,-t))\\[.1in]
    A_o(x,t) & =  \sum_{k\geq 0} F_{2k+1}(x)t^{2k+1}\\
      & =  \frac 12(A(x,t)-A(x,-t)). \end{split}\end{equation}
Multiply equation (\ref{eq:akx}) by $t^k$ and sum on $k\geq 0$. We
obtain
  \beq \frac{\partial A(x,t)}{\partial x} = tA_e(x,t)A(x,t)+
    A_o(x,t)A(x,t). \label{eq:paxt} \eeq
Substituting $-t$ for $t$ yields
    \beq \frac{\partial A(x,-t)}{\partial x} = -tA_e(x,t)A(x,-t)-
    A_o(x,t)A(x,-t). \label{eq:paxt2} \eeq
Adding and subtracting equations (\ref{eq:paxt}) and (\ref{eq:paxt2})
gives the following system of differential equations for
$A_e=A_e(x,t)$ and $A_o=A_o(x,t)$:
  \begin{align} \frac{\partial A_e}{\partial x} & = tA_eA_o+A_o^2
     \label{eq:diff1}\\[.05in] \frac{\partial A_o}{\partial x} & =
       tA_e^2+A_eA_o. \label{eq:diff2} \end{align}
Thus we need to solve this system of equations in order to find
$A(x,t) = A_e(x,t)+A_o(x,t)$.

\begin{theorem} \label{thm:main}
We have
      \begin{align} B(x,t) & = \frac{1+\rho+2te^{\rho x}+(1-\rho)
     e^{2\rho x}}{1+\rho-t^2+(1-\rho-t^2)e^{2\rho x}}
    \label{eq:fxt2}\\[.05in] A(x,t) & = (1-t)B(x,t)
      \label{eq:fxt3}\\
     & = (1-t)\frac{1+\rho+2te^{\rho x}+(1-\rho)
     e^{2\rho x}}{1+\rho-t^2+(1-\rho-t^2)e^{2\rho x}},
      \label{eq:fxt4} \end{align}
where $\rho=\sqrt{1-t^2}$. 
\end{theorem}

\proof
We can simply verify that the stated expression (\ref{eq:fxt4}) for
$A(x,t)$ satisfies (\ref{eq:diff1}) and (\ref{eq:diff2}) with the
initial condition $A(0,t)=1$, a routine computation (especially with
the use of a computer). The relationship (\ref{eq:fxt3}) between
$A(x,t)$ and $B(x,t)$ is then an immediate consequence of
(\ref{eq:bkak}), which is equivalent to $a_k(n)=b_k(n)-b_k(n-1)$.

It might be of interest, however, to explain how the formula
(\ref{eq:fxt4}) for $A(x,t)$ can be derived if the answer is not known
in advance. If we divide equation (\ref{eq:diff1}) by
(\ref{eq:diff2}), then we obtain
  \[ \frac{\partial A_e/\partial x}{\partial A_o/\partial x} =
     \frac{A_o}{A_e}. \]
Hence $\frac{\partial}{\partial x}(A_e^2-A_o^2)=0$, so $A_e^2-A_o^2$
is independent of $x$. This observation suggests computing the
generating function in $t$ for $A_e^2-A_o^2$, which the computer
shows is equal to $1+O(t^N)$ for a large value of $N$. Assuming then
that $A_e^2-A_o^2=1$ (or even proving it combinatorially), we can
substitute $\sqrt{1-A_e^2}$ for $A_o$ in (\ref{eq:diff1}) to obtain
  \[ \frac{\partial A_e}{\partial x} = tA_e\sqrt{A_e^2-1}+A_e^2-1, \]
a single differential equation for $A_e$. This equation can routinely
be solved by separation of variables (though some care must be taken
to choose the correct branch of the resulting integral, including the
correct sign of $\sqrt{A_e^2-1}$); we will spare the reader the
details. A similar argument yields $A_o$, so we obtain $A=A_e+A_o$.
\qed

\medskip
\textsc{Note.} Ira Gessel has pointed out the following simplified
expression for $B(x,t$):
  \beq B(x,t) = \frac{ 2/\rho}{\displaystyle
      1-\frac{1-\rho}t e^{\rho 
      x}} -\frac1{\sqrt{1-t^2}}. \label{eq:bgessel} \eeq

\section{Consequences.} \label{sec3}
A number of corollaries follow from Theorem~\ref{thm:main}. The first
is the explicit expressions for $a_k(n)$ and $b_k(n)$ stated in the
introduction. I am grateful to Ira Gessel for providing the proof
given below.

\begin{corollary} \label{cor:aknbkn}
For all $k,n\geq 1$ we have
  \begin{align} 
  b_k(n) & = \frac{1}{2^{k-1}} \sum_{\substack{r+2s\leq k\\ r\equiv
     k\,(\mathrm{mod}\,2)}} 
     (-2)^s{k-s\choose (k+r)/2}{n\choose s}r^n \label{eq:bkns}\\
   a_k(n) & = b_k(n)-b_{k-1}(n). \label{eq:akns} \end{align}
\end{corollary}

\proof
Define $b'_k(n)$ to be the right-hand side of \eqref{eq:bkns}, and
set
  \[ B'(x,t) = \sum_{k,n\geq 0} b'_k(n)t^k\frac{x^n}{n!}. \]
Set $n=s+m$ and $k=r+2s+2l$, so 
\begin{align*}B'(x,t) &= \sum_{r,s,l,m} (-1)^s 2^{1-r-s-2l}
  \binom{r+s+2l}{r+s+l} \binom{s+m}s r^{s+m} t^{r+2s+2l} 
\frac {x^{s+m}}{(s+m)!}\\
&=2\sum_{r,s\ge0} \left(\frac t2\right)^r \frac{(-rt^2x/2)^s}{s!}
\biggl[\sum_l \binom{r+s+2l}l \left(\frac{t^2}4\right)^l\biggr]
\biggl[\sum_m \frac{(rx)^m}{m!}\biggr].
\end{align*}
The sum on $m$ is $e^{rx}$. Using the formula
\[\sum_k \binom{2k+a}{k}u^k = \frac{C(u)^a}{\sqrt{1-4u}},\]
where 
 \[ C(u) = \sum_{n\geq 0}C_n u^n =\frac{1-\sqrt{1-4x}}{2x}, 
 \] 
the generating function for the Catalan numbers
$C_n=\frac{1}{n+1}\binom{2n}{n}$, we find that the sum on  $l$ is 
\[\frac {C(t^2/4)^{r+s}}{\sqrt{1-t^2}} = \frac 1\rho \left(\frac
  {2-2\rho}{t^2}\right)^{r+s}.\] 
 
Thus 
\begin{align*}
B'(x,t) &= \frac 2\rho 
\sum_{r,s\ge0} \left(\frac t2\right)^r \frac{(-rt^2x/2)^s}{s!}
e^{rx}\left(\frac {2-2\rho}{t^2}\right)^{r+s}\\
 & =  \frac 2\rho \sum_{r}\left(\frac {1-\rho} t e^x\right)^r 
   \sum_s \frac{(-r(1-\rho)x)^s}{s!}\\
   &= \frac 2\rho \sum_{r}\left(\frac {1-\rho} t e^x\right)^r 
        e^{-r(1-\rho)x}\\
   &= \frac 2\rho \frac1 {\displaystyle 1-\frac{1-\rho}t e^{\rho x}},
\end{align*}
and the proof of \eqref{eq:bkns} follows from
(\ref{eq:bgessel}). Equation~\eqref{eq:akns} is then an immediate
consequence of \eqref{eq:bkak}.
\qed

\medskip
By Corollary~\ref{cor:aknbkn}, when $k$ is fixed $b_k(n)$
is a linear combination of $k^n$, $(k-2)^n$, $(k-4)^n, \dots$ with
coefficients that are polynomials in $n$. For $k\leq 6$ we have 
  \begin{align*}  b_2(n) & = 2^{n-1}\\
        b_3(n) & = \frac 14(3^n-2n+3)\\
    b_4(n) & = \frac 18(4^n-2(n-2)2^n)\\ 
    b_5(n) & = \frac{1}{16}(5^n-(2n-5)3^n+2(n^2-5n+5)) \\
    b_6(n) & = \frac{1}{32}(6^n-2(n-3)4^n+(2n^2-12n+15)2^n). 
\end{align*}

As a further application of Theorem~\ref{thm:main} we can obtain the
factorial moment generating function
  \[ F(x,t) = \sum_{s,n\geq 0} \nu_j(n)x^n\frac{t^j}{j!}, \]
where
  \[ \nu_j(n) = \frac{1}{n!}\sum_{w\in\sn}(\as(w))_j 
     = \frac{1}{n!}\sum_k a_k(n)(k)_j. \]
and 
  \[ (h)_j =h(h-1)\cdots (h-j+1). \]
Namely, we have
  \begin{align*} \left.\frac{\partial^j A(x,t)}{\partial t^j}\right|_{t=1}
    & = \sum_{n\geq 0}\frac{1}{n!}\sum_{k\geq 0}a_k(n)(k)_j\,x^n\\
   & = \sum_{n\geq 0}\nu_j(n)x^n. \end{align*}
On the other hand, by Taylor's theorem we have
  \[ A(x,t) = \sum_{j\geq 0} \left.\frac{\partial^j A(x,t)}{\partial
     t^j}\right|_{t=1} \frac{(t-1)^j}{j!}. \]
It follows that \
  \beq F(x,t)=A(x,t+1). \label{eq:fa} \eeq
(Note that it is not at all a priori obvious from the form of
$A(x,t+1)$ obtained by substituting $t+1$ for $t$ in \eqref{eq:fxt4}
that it even has a Taylor series expansion at $t=0$.) From
equations (\ref{eq:fxt4}) and (\ref{eq:fa}) it is easy to compute
(using a computer) the generating functions
  \[ M_j(x) =\sum_{n\geq 0} \nu_j(n)x^n \]
for small $j$. For $1\leq j\leq 4$ we get

 \begin{align*} M_1(x) & = \frac{6x-3x^2+x^3}{6(1-x)^2}\\[.05in]
    M_2(x) & = \frac{90x^2-15x^4+6x^5-x^6}{90(1-x)^3}\\[.05in]
    M_3(x) & = \frac{2520x^3-315x^4+189x^5-231x^6+93x^7-18x^8+2x^9}
                    {1260(1-x)^4}\\[.05in]
   M_4(x) & =
    \frac{N_4(x)}{9450(1-x)^5}, \end{align*}
where
  \[ N_4(x)=47250x^4-3780x^6+2880x^7-2385x^8+1060x^9-258x^{10} \]
\vspace{-.3in}
  \[ \hspace{-1.5in}  +36x^{11}-3x^{12}. \]
It is not difficult to see that in general $M_j(x)$ is a rational
function of $x$ with denominator $(1-x)^{j+1}$. It follows from
standard properties of rational generating functions
\cite[{\S}4.3]{ec1} that for fixed $j$ we have that $\nu_j(n)$ is a
polynomial in $n$ of degree $j$ for $n$ sufficiently large. In
particular, we have
  \begin{align} \nu_1(n) & = \frac{4n+1}{6},\ n\geq 2
       \label{eq:mean}\\ 
    \nu_2(n) & = \frac{40n^2-24n-19}{90},\ n\geq 4\nonumber \\
    \nu_3(n) & =
    \frac{1120n^3-2856n^2+440n+1581}{3780},\ n\geq
    6. \nonumber \end{align}
Note in particular that $\nu_1(n)$ is just the expectation (mean) of
$\as_n$. The simple formula $(4n+1)/6$ for this quantity should be
contrasted with the situation for the length $\is_n(w)$ of the longest
increasing subsequence of $w\in\sn$, where even the asymptotic
formula $E(n)\sim 2\sqrt{n}$ for the expectation is a highly
nontrivial result \cite[{\S}3]{ids}. A simple proof of \eqref{eq:mean}
follows from \eqref{eq:tg} and an argument of Knuth
\cite[Exer.~5.1.3.15]{knuth}. 

From the formulas for $\nu_1(n)$ and $\nu_2(n)$ we easily compute the
variance var$(\as_n)$ of $\as_n$, namely,
  \beq \mathrm{var}(\as_n) = \nu_2(n)+\nu_1(n)-\nu_1(n)^2 = 
    \frac{32n-13}{180},\ n\geq 4. \label{eq:var} \eeq
\indent We now consider a further application of
Theorem~\ref{thm:main}. Let 
  \beq T_n(t) = \sum_{k=0}^n a_k(n)t^k. \label{eq:tnt} \eeq
For instance,
  \begin{align*} T_1(t) & = t\\
   T_2(t) & = t+t^2\\
   T_3(t) & = t+3t^2+2t^3\\
   T_4(t) & = t+7t^2+11t^3+5t^4\\
   T_5(t) & = t+15t^2+43t^3+45t^4+16t^5\\
   T_6(t) & = t+31t^2+148t^3+268t^4+211t^5+61t^6\\
   T_7(t) & = t+63t^2+480t^3+1344t^4+1767t^5+1113t^6+272t^7. \end{align*}

\begin{corollary} \label{cor:div}
The polynomial $T_n(t)$ is divisible by $(1+t)^{\lfloor
n/2\rfloor}$. Moreover, if $U_n(t)=T_n(t)/(1+t)^{\lfloor
n/2\rfloor}$, then
  \[ U_{2n}(-1) = -U_{2n+1}(-1)=\frac{(-1)^n E_{2n+1}}{2^n}, 
  \]
where $E_{2n+1}$ denotes a tangent number.
\end{corollary}

\proof
Let $A_e(x,t)$ and $A_o(x,t)$ be the even and odd parts of $A(x,t)$ as
in equation~\eqref{eq:aeao}. By the definition of $A_e(x)$ we have
  \[ A_e(x/\sqrt{1+t},t) = \sum_{n\geq
    0}\frac{T_{2n}(t)}{(1+t)^n}\frac{x^{2n}}{(2n)!}. \] 
With the help of the computer we compute that
  \begin{align*} \lim_{t\to -1} A_e(x/\sqrt{1+t},t) & = \mathrm{sech}^2
    \frac{x}{\sqrt{2}}\\ & = \sum_{n\geq 0}\frac{(-1)^n E_{2n+1}}
           {2^n}\frac{x^{2n}}{(2n)!}. \end{align*}
Hence the desired result is true for $T_{2n}(t)$. Similarly,
   \begin{align*} \lim_{t\to -1} \sqrt{1+t}\,A_o(x/\sqrt{1+t},t) & =
   -\sqrt{2}\tanh \frac{x}{\sqrt{2}}\\ & = -\sum_{n\geq
     0}\frac{(-1)^n E_{2n+1}} {2^n}\frac{x^{2n+1}}{(2n+1)!}, \end{align*}
proving the result for $T_{2n+1}(t)$.
\qed

\medskip
By Corollary~\ref{cor:div} we have $T_n(-1)=0$ for $n\geq 2$. In other
words, for $n\geq 2$ we have
  \[ \#\{w\in\sn\st \as_n(w)\ \mathrm{even}\} =
     \#\{w\in\sn\st \as_n(w)\ \mathrm{odd}\} =
     \frac{n!}{2}. \]
A simple combinatorial proof of this fact follows from switching the
last two elements of $w$; it is easy to see that this operation either
increases or decreases $as_n(w)$ by 1, as first pointed out by
M. B\'ona and P. Pylyavskyy. More generally, a combinatorial proof of
Corollary~\eqref{cor:div} is a consequence of equation~\eqref{eq:tg}
below and an argument of B\'ona \cite[Lemma~1.40]{bona}.

The formulas \eqref{eq:mean} and \eqref{eq:var} for the mean and
variance of $\as_n$ suggest in analogy with \eqref{eq:bdj} that
$\as_n$ will have a limiting distribution $K(t)$ defined by 
  \[ K(t) = \lim_{n\to\infty} \mathrm{Prob}\left(
   \frac{\as_n(w)-2n/3}{\sqrt{n}} \leq t\right), \]
for all $t\in\rr$, where $w$ is chosen uniformly from $\sn$. Indeed,
we have that  $K(t)$ is a Gaussian distribution with variance
$8/45$: 
    \beq K(t) = \frac{1}{\sqrt{\pi}} \int_{-\infty}^{t\sqrt{45}/4}
        e^{-s^2}ds. \label{eq:kt} \eeq
It was pointed out by Pemantle (private communication) that
equation~\eqref{eq:kt} is a consequence of the result \cite[Thms.~3.1,
3.3, or 3.5]{p-w} and possibly also \cite{b-r}. An independent proof
was also given by Widom \cite{widom}, and in the next section we
explain an additional method of proof.

\section{Relationship to alternating runs.} \label{sec4}
A \emph{run} of a permutation $w=a_1\cdots a_n\in\sn$ is a maximal
factor (subsequence of consecutive elements) which is increasing. An
\emph{alternating run} is a maximal factor that is increasing or
decreasing. (Perhaps ``birun'' would be a better term.) For instance,
the permutation 64283157 has four alternating runs, viz., 642, 28,
831, and 157. Let $g_k(n)$ be the number of permutations $w\in\sn$
with $k$ alternating runs. It is easy to see, as pointed out by B\'ona
\cite{bona2}, that
  \beq a_k(n) = \frac 12(g_{k-1}(n)+g_k(n)),\ \ n\geq 2. 
    \label{eq:ag} \eeq
If we define $G_n(t)=\sum_k g_k(n)t^k$, then equation~\eqref{eq:ag} is
equivalent to the formula
  \beq T_n(t) = \frac 12(1+t)G_n(t), \label{eq:tg} \eeq
where $T_n(t)$ is defined by \eqref{eq:tnt}. 

Research on the numbers $g_k(n)$ go back to the nineteenth
century; for references see Bona \cite[{\S}1.2]{bona} and Knuth
\cite[Exer.~5.1.3.15--16]{knuth}. In particular, let $A_n(t)$ denote
the $n$th \emph{Eulerian polynomial}, i.e., 
  $$ A_n(t) = \sum_{w\in\sn} t^{1+\mathrm{des}(w)}, $$
where des$(w)$ denotes the number of descents of $w$ (the size of the
descent set defined in equation~(\ref{eq:deset})). It was shown by
David and Barton \cite[pp.~157--162]{d-b} and stated more concisely by
Knuth \cite[p.~605]{knuth} that
  \[ G_n(t) = \left( \frac{1+t}{2}\right)^{n-1}
    (1+w)^{n+1}A_n\left( \frac{1-w}{1+w}\right),\ \ n\geq 2, \]
where $w =\sqrt{\frac{1-t}{1+t}}$. Theorem~\ref{thm:main} is then a
straightforward consequence of the well-known generating function
(e.g., \cite[Thm.~1.7]{bona})
  \[ \sum_{n\geq 0} A_n(t)\frac{x^n}{n!} =
      \frac{1-t}{1-te^{(1-t)x}}. \] 
\indent It is also well-known (e.g., \cite[Thm.~1.10]{bona}) that the
Eulerian polynomial $A_n(t)$ has only real zeros, and that the zeros
of $A_n(t)$ and $A_{n+1}(t)$ interlace. From this fact Wilf
\cite{wilf} showed that the polynomials $G_n(t)$ have (interlacing)
real zeros, and hence by \eqref{eq:tg} the polynomials $T_n(t)$ also
have real zeros. It is then a consequence of standard results (e.g.,
\cite[Thm.~2]{bender}) that the numbers $a_k(n)$ for fixed $n$ are
asymptotically normal as $n\rightarrow \infty$, yielding another proof
of \eqref{eq:kt}.

\section{Open problems.}
In this section we mention three directions of possible generalization
of our work above.
\begin{enumerate}
\item Let $\is(m,w)$ denote the length of the longest subsequence of
  $w\in\sn$ that is a union of $m$ increasing subsequences, so
  $\is(w)=\is(1,w)$. The numbers $\is(m,w)$ have many interesting
  properties, summarized in \cite[{\S}4]{ids}. Can anything be said
  about the analogue for alternating sequences, i.e., the length
  $\as(m,w)$ of the longest subsequence of $w$ that is a union of $m$
  alternating subsequences?  This question can also be formulated in
  terms of the lengths of the alternating runs of $w$.
 \item Can the results for increasing subsequences and alternating
 subsequences be generalized to other ``patterns''? More specifically,
 let $\sigma$ be a (finite) word in the letters $U$ and $D$, e.g.,
 $\sigma=UUDUD$. Let $\sigma^\infty$ denote the infinite word
 $\sigma\sigma\sigma\cdots$, e.g., 
   \[ (UUD)^\infty = UUDUUDUUD\cdots. \]
For this example, we have for instance that $UUDUUDU$ is a prefix of
$\sigma^\infty$ of length 7. 

Let $\tau=a_1 a_2 \cdots a_{m-1}$ be a word of length $m-1$ in the
letters $U$ and $D$. A sequence $v=v_1 v_2\cdots v_m$ of integers is
said to have \emph{descent word} $\tau$ if $v_i>v_{i+1}$ whenever
$a_i=D$, and $v_i<v_{i+1}$ whenever $a_i=U$. Thus $v$ is increasing if
and only if $\tau=U^{m-1}$, and $v$ is alternating if and only if
$\tau= (DU)^{j-1}$ or $\tau=(DU)^{j-1}D$ depending on whether $m=2j-1$
or $m=2j$.

Now let $w\in\sn$ and define $\lgn_\sigma(w)$ to be the length of
longest subsequence of $w$ whose descent word is a prefix of
$\sigma^\infty$. Thus $\lgn_U(w)=\is_n(w)$ and $\lgn_{DU}(w) =
\as_n(w)$. What can be said in general about $\lgn_\sigma(w)$? In
particular, let
  \[ E_\sigma(n) = \frac{1}{n!}\sum_{w\in\sn} \lgn_\sigma(w), \]
the expectation of $\lgn_\sigma(w)$ for $w\in\sn$. Note that
$E_U(n)\sim 2\sqrt{n}$ by \eqref{eq:en}, and
$E_{DU}(n)\sim 2n/3$ by \eqref{eq:mvar}. Is it true that for any
$\sigma$ we have $E_\sigma(n)\sim \alpha n^c$ for some $\alpha,c>0$?
Or at least that for some $c>0$ (depending on $\sigma$) we have
  \[ \lim_{n\to\infty} \frac{\log E_\sigma(n)}{\log n} =c, \]
in which case can we determine $c$ explicitly?
\item The \emph{descent set} $D(w)$ of a permutation $w=w_1\cdots w_n$
  is defined by
  \beq D(w) = \{ i\st w_i>w_{i+1}\} \subseteq [n-1], \label{eq:deset}
    \eeq
where $[n-1]=\{1,2,\dots,n-1\}$. Thus $w$ is alternating if and only
if $D(w)=\{1,3,5,\dots\}\cap [n-1]$. Let $S\subseteq [k-1]$. What can
be said about the number $b_{k,S}(n)$ of permutations $w\in \sn$ that
avoid all $v\in\fs_k$ satisfying $D(v)=S$? In particular, what is the
value $L_{k,S}=\lim_{n\to\infty} b_{k,S}(n)^{1/n}$? (It follows from
\cite{arratia} and \cite{m-t}, generalized in an obvious way, that
this limit exists and is finite.) For instance, if $S=\emptyset$  
or $S=[k-1]$, then it follows from \cite{regev} that
$L_{k,S}=(k-1)^2$. On the other hand, if $S=\{1,3,5,\dots\}\cap [k-1]$
then it follows from \eqref{eq:bkns} that $L_{k,S}=k-1$.
\end{enumerate}

\pagebreak

\end{document}